\documentclass[a4paper,12pt]{article}
\usepackage[latin1]{inputenc}
\usepackage{amsmath}
\usepackage{amssymb}
\usepackage{amscd}
\usepackage{array}
\usepackage[all]{xy}
\usepackage{a4}
\DeclareMathAlphabet{\mathfr}{U}{euf}{m}{n}
\newtheorem{teo}{Theorem}[section]

\newtheorem{rem}{Remark}[section]

\newcommand{\Q}{\mathbb Q}
\newcommand{\Qbar}{{\overline{\mathbb Q}}}
\newcommand{\Gal}{\mathrm{Gal}}

\newcommand{\Z}{\mathbb Z}

\newcommand{\F}{\mathbb F}

\newcommand{\PGL}{\mathrm{PGL}}
\newcommand{\PSL}{\mathrm{PSL}}
\newcommand{\eps}{\varepsilon}
\newcommand{\To}{\longrightarrow}

\newcommand{\Aut}{\operatorname{Aut}}

\newcommand{\G}{\mathrm{G}}

\title{Plane quartic twists of $X(5,3)$}
\author{Julio Fern\'andez, Josep Gonz\'alez, Joan-C. Lario}
\date{\today}
\begin{document}
\maketitle


\begin{abstract}
Given an odd surjective Galois representation $\varrho\,\colon\!\G_\Q\to\PGL_2(\F_3)$ and a
positive integer $N$\!, there exists a twisted modular curve $X(N,3)_\varrho$
defined over $\Q$ whose rational points classify the quadratic $\Q$-curves of degree $N$
realizing~$\varrho$.
The paper gives a method to provide an explicit plane quartic model for this curve
in the genus-three case $N\!=5$.
\end{abstract}


\section{Introduction}

Let $p$ be an odd prime. The so-called $\Q$-curves (non-CM elliptic curves over $\Qbar$ that are isogenous
to all their Galois conjugates) are a source of representations of the absolute Galois group $\G_\Q$
into $\PGL_2(\F_p)$. We refer to \cite{Fe-La-Rio} or \cite{Fe} for a detailed construction of the odd
representation $\varrho_{E,p}\,\colon\G_\Q\to\PGL_2(\F_p)$ arising from the $p$-torsion of a $\Q$-curve~$E$.
It seems natural to ask for the {\it frequency} of such representations $\varrho_{E,p}$\, among all odd
$2$-dimensional projective mod $p$\, Galois representations. We say that a \mbox{$\Q$-curve $E$} 
\emph{realizes} a given representation $\varrho\,\colon\G_\Q\to\PGL_2(\F_p)$\, if
\,$\varrho_{E,p}=\varrho$\, up to conjugation in $\PGL_2(\F_p)$.

The moduli problem classifying the quadratic $\Q$-curves that
realize such a representation $\varrho$\, is the subject
of the PhD thesis \cite{julio}. The problem
splits into different cases depending on the value mod $p$\,
of the eventual degrees $N$ for the isogeny from the $\Q$-curve
to its Galois conjugate. For a fixed $N$\!, the quadratic $\Q$-curves
of degree~$N$ realizing $\varrho$ are given by the
rational points on some twist of a modular curve. The aim
of this paper is to explain a method to obtain good models
over $\Q$ for some of these twisted curves, namely for those
corresponding to the octahedral genus-three case $N\!=5$, $p=3$.

In Section~2 we briefly introduce the modular curve $X(N,p)$ and its twists $X(N,p)_\varrho$.
We refer to \cite{Fe} for this section.
In Section~3 we give a plane quartic model for $X(5,3)$ along with an explicit description of the
natural map $X(5,3)\rightarrow X^+(5)$. An algorithm to produce a plane quartic model
for $X(5,3)_\varrho$\, whenever $\varrho$\, is surjective is then explained in Section~4, where the input
is given by a degree-four polynomial in $\Z[X]$ having the same splitting field as $\varrho$.
Finally, in Section~5 we present an example illustrating that all computations can be carried out
by any of the available standard algebraic manipulation systems. We also compute \emph{by hand} some
rational points on the resulting quartic model and show that the Chabauty-Coleman method fails
to be of use in this case. The authors expect that some other developing techniques to
determine the set of all rational points on small genus curves
can be applied to the models provided by the method described in the paper.


\section{The twisted modular curve $X(N,p)_\varrho$}

Let $N\!>1$ be an integer prime to $p$. We denote by~$X(N,p)$ the fiber product over~$X(1)$ of 
the modular curves $X_0(N)$ and $X(p)$. We take for $X_0(N)$ its canonical model over $\Q$.
As for $X(p)$, we fix the rational model attached to a matrix~$V$ in
$\PGL_2(\F_p)\setminus\PSL_2(\F_p)$ of order $2$, as a particular case of a general procedure that 
can be found in Section II.3 of \cite{Li} and Section~2 of \cite{Ma77}. Its $\Q$-isomorphism class 
does not depend on the choice of such a matrix.
We denote by $\mathcal{W}(N,p)$ the automorphism group of the covering
$$X(N,p)\To X^+(N),$$
where $X^+(N)$ is the quotient of $X_0(N)$ by the Atkin-Lehner involution $w_{\scriptscriptstyle N}$.
We recall that the 
rational points on~$X^+(N)$ yield the isomorphism classes
of quadratic $\Q$-curves of degree~$N$ up to Galois conjugation.

Assume $N$ to be a non-square mod $p$. The automorphism group~$\mathcal{W}(N,p)$
is then canonically isomorphic to $\PGL_2(\F_p)$. If we put $w$ for the involution on~$X(N,p)$
corresponding through this isomorphism to the above matrix~$V$ and identify $\mathcal{W}(N,p)$
with its (inner) automorphism group, the Galois action on~$\mathcal{W}(N,p)$ is given by the
morphism
$$\eps\,\colon\G_\Q\,\To\,\F_p^{\,*}/\,{\F_p^{\,*}}^2\,
\simeq\,\langle w\rangle\,\hookrightarrow\,\mathcal{W}(N,p)$$
obtained from the mod $p$\, cyclotomic character $\G_\Q\!\To\F_p^{\,*}$.

Suppose that we are now given a surjective Galois representation
$$\varrho\,\colon\,\G_\Q\,\To\,\PGL_2(\F_p)$$
with non-cyclotomic determinant. A $\Q$-curve of degree $N$
realizing $\varrho$ must then be defined over the quadratic field attached to
the Galois character
$$\eps\,\det\varrho\,\colon\G_\Q\,\To\,\F_p^{\,*}/\,{\F_p^{\,*}}^2\,\simeq\,\{\pm 1\}\,.$$
For the moduli classification of such \mbox{$\Q$-curves,} we produce a twist of $X(N,p)$ 
from a certain element in the cohomology set $H^1\!\left(\G_\Q, \mathcal{W}(N,p)\right)$. 
Specifically, we take the $1$-cocycle \,$\xi=\varrho\,\eps$, where we view $\varrho$\, as 
a map onto $\mathcal{W}(N,p)$ through the above canonical isomorphism:
$$\varrho\,\colon\G_\Q\,\To\,\PGL_2(\F_p)\,\stackrel{\simeq}\To\,\mathcal{W}(N,p)\,.$$
For the twist of~$X(N,p)$ attached to~$\xi$, we fix a rational model $X(N,p)_\varrho$\, along with an
isomorphism
$$\Psi\,\colon\,X(N,p)_\varrho\,\To\,X(N,p)$$
satisfying \,$\Psi\,=\,\xi_\sigma{^\sigma\!}\Psi$\, for every~$\sigma$ in~$\G_\Q$.
We note that the $\Q$-isomorphism class of $X(N,p)_\varrho$\, is an
invariant of the conjugacy class of $\varrho$. In other words, it only depends on the splitting
field of $\varrho$, since every automorphism of $\PGL_2(\F_p)$ is inner.

\begin{teo}\label{non-cyclo} There exists a quadratic \,$\Q$-curve of degree $N$ realizing $\varrho$ if and only
if the set of non-cuspidal non-CM \,rational points on the curve $X(N,p)_\varrho$ is nonempty. In this case,
the composition
$$X(N,p)_\varrho \stackrel{\Psi}\To X(N,p) \To X^+(N)$$
defines a one-to-one correspondence between this set of points and the set of
isomorphism classes of quadratic \,\mbox{$\Q$-curves} of degree~$N$ up to Galois conjugation
realizing~$\varrho$.
\end{teo}


\begin{rem}
The genus of $X(N,p)$ is never two, and the only genus-three case appears for $N\!=5$ and $p=3$.
\end{rem}


\section{A plane quartic model for $X(5,3)$}
The automorphism of the complex upper-half plane given by \,$\tau \mapsto \tau/3$\,
induces an isomorphism
$$\Phi\,\colon X(N,3) \To X_0(9N)$$\\[-10pt]
defined over $\Q$, so that $\Phi^*\!\left(\,\Q\left(X_0(9N)\right)\right)$ is the function field of the rational
model for $X(N,3)$ fixed in the previous section.
Moreover, the above involution $w$ on $X(N,3)$ can be chosen to correspond through $\Phi$ to the 
Atkin-Lehner involution $w_{\scriptscriptstyle N}$ on
$X_0(9N)$. From now on, we take~$N=5$ and ease the notation as follows: for a
function~$x\in\Q\left(X_0(45)\right)$, a regular differential $\omega\in\Omega^1\!\left(X_0(45)\right)$ or an
automorphism $W\in\Aut\left(X_0(45)\right)$, we put $\overline x=\Phi^*(x)$, $\overline\omega=\Phi^*(\omega)$,
$\overline W=\Phi^{-1}\,W\,\Phi$ for the corresponding function, differential or automorphism on $X(5,3)$.


\subsection{An  equation for  $X_0(45)$}
The jacobian of $X_0(45)$  is $\Q$-isogenous to \,$J_0(15)^2\times J_0(45)^{\operatorname{new}}$\!, 
where $J_0(15)$ and~$J_0(45)^{\operatorname{new}}$ are  elliptic curves over $\Q$ of conductors $15$ and $45$,
respectively. A basis for the $\Q$-vector space \,$\Omega^1_\Q\!\left(X_0(45)\right)$ is
given by
$$\omega_1=f_1(q)\,\frac{\,dq\,}{\,q\,}\,, \ \ \ \ \omega_2=f_1(q^3)\,\frac{\,dq\,}{\,q\,}\,, \ \ \ \ \omega_3=f_2(q)\,\frac{\,dq\,}{\,q\,}\,,
$$
where $f_1$ and $f_2$ are the normalized weight-two newforms of levels $15$ and $45$, respectively:\\[-10pt]
$$f_1 = q - q^2 - q^3 - q^4 + q^5 + q^6 + 3q^8 + q^9 - q^{10} - 4q^{11} + q^{12} - 2q^{13} - q^{15} - q^{16} +
2q^{17}+\cdots
$$
\vspace{-6truemm}
$$
f_2 = q + q^2 - q^4 - q^5 - 3q^8 - q^{10} + 4q^{11} - 2q^{13} - q^{16} - 2q^{17}+\cdots\ \ \
\ \ \ \ \ \ \ \ \ \ \ \ \ \ \ \ \ \
$$

The genus-three curve  $X_0(45)$ is nonhyperelliptic (see \cite{Ogg}), so the image of $X_0(45)$ under the
canonical embedding is the zero locus of a homogenous polynomial $P$ in $\Q [X,Y,Z]$ of degree $4$. Such a
polynomial $P$\!, unique up to non-zero rational multiples, satisfies \,$P(\omega_1,\omega_2,\omega_3)=0$\, 
and can be explicitly determined using the first seventeen Fourier coefficients of each $\omega_i$\, 
(see Section~2 of \cite{bggp}).
This yields the following affine equation for $X_0(45)$\,:
\begin{equation}\label{X045}
x^4 - 2\,x^2\,y^2 + 81\,y^4 - 2\,x^2 - 16\,x\,y - 18\,y^2 + 1=0,
\end{equation}
where \,$x=\omega_1/\omega_3$\, and \,$y=\omega_2/\omega_3$. 
In particular, \,$\Q\left(X(5,3)\right)=\Q (\overline x,\overline y)$.


\subsection{The group  $\Aut\left(X_0(45)\right)$}
The group $\Aut\left(X_0(45)\right)$ is generated by the Atkin-Lehner
involutions \,$w_5, w_9$\, and the automorphism \,$S$\, induced by the map \,$\tau\mapsto \tau+1/3$\,
on the complex upper-half plane (see \cite{Ke-Mo}, \cite{Le-Ne}). The action of these generators 
on the $\Qbar$-vector space $\Omega^1\!\left(X_0(45)\right)$ is displayed in the following table:
$$
\begin{array}{c|c|c|c|}
   & w_5 & w_9& S\\[2pt]
   \hline
   \hline
   & & & \\[-7pt]
\ \omega_1 \ & \ -\omega_1 \ & \ \,\,\,\,\,3\,\omega_2  \ & 
\ -1/2\,\omega_1-3/2\,\omega_2+\sqrt{-3\,}/2\,\omega_3 \ \\[7pt]
\ \omega_2 \ & \ -\omega_2 \  &  \ 1/3\,\omega_1 \ &  \ \omega_2 \ \\[7pt]
\ \omega_3 \ & \ \,\,\,\,\,\omega_3 \ & \ \,\,\,\,\,-\omega_3 \ & \ \sqrt{-3\,}/2\,\omega_1+\sqrt{-3\,}/2\,\omega_2-1/2\,\omega_3\ \\[5pt]
\hline
\end{array}
$$
It can be easily checked from these relations that $\Aut\left(X_0(45)\right)$ has the same order as 
$\PGL_2(\F_3)$, so that \,$\Aut\left(X(5,3)\right) = \mathcal{W}(5,3) = 
\langle\,\overline w_5,\,\overline w_9,\,\overline S\,\rangle \simeq \PGL_2(\F_3)$.
Let us also note that $\Aut_{\Q}\left(X_0(45)\right)$ is the subgroup of Atkin-Lehner
involutions, while the automorphism $S$ is defined over the quadratic field $\Q(\sqrt{-3\,}\,)$ 
and satisfies \,${}^{\nu\!}S=S^2$\, for the non-trivial automorphism $\nu\in\Gal(\Q(\sqrt{-3\,}\,)/\Q)$. 
In particular, the group $\langle w_9\,S\,w_9\rangle$ is $\G_{\Q}$-stable.
As a matter of fact, the degree-three covering $X_0(45)\To X_0(45)/\langle w_9\,S\,w_9\rangle$
is the natural projection $X_0(45)\To X_0(15)$, since the pullback of  
\,$\Omega^1(X_0(45)/\langle w_9\,S\,w_9\rangle)$\, is \,$\langle\omega_1\rangle$.

\subsection{The covering $X(5,3)\rightarrow X^+(5)$}
Consider the composition \,$X_0(45)\,\stackrel{\pi_1}\To\,X_0(15)\,\stackrel{\pi_2}\To\,X_0(5)\,\stackrel{\pi_3}\To\,X^+(5)$,
where $\pi_1$,  $\pi_2$ and $\pi_3$ are the natural projections. We fix for $X_0(45)$ the model given by
equation (\ref{X045}). For~$X_0(15)$, we take the minimal equation
\begin{equation}\label{minimal}
v^2+u\,v+v=u^3+u^2-10\,u-10
\end{equation}
given in \cite{cremona97}, where the functions $u,v\in\Q\left(X_0(15)\right)$ have a unique pole at the cusp
$\infty$ with respective multiplicities $2$ and $3$\,:
$$
u=\frac{1}{\,q^2\,}+\frac{1}{\,q\,}+1 + 2\, q + 4\, q^2  + \cdots\, \ \ \ \
v=\frac{1}{\,q^3\,}+\frac{1}{\,q^2\,}+\frac{2}{\,q\,}+3 + 2\, q + 5\, q^2 +\cdots
$$
As for the function fields of $X_0(5)$ and $X^+(5)$, we take the following generators over~$\Q$,
respectively:
$$
\begin{array}{rcl}
G(\tau) = & {\displaystyle \left(\frac{\eta (\tau)}{\,\eta (5\tau)\,}\right)^{\!6}} & = \
{\displaystyle \frac{1}{\,q\,} - 6 + 9\, q + 10\, q^2 +\cdots}\\[15pt]
t(\tau) = & {\displaystyle G(\tau)+\frac{5^3}{~\,G(\tau)\,}} & = \ 
{\displaystyle \frac{1}{\,q\,} - 6 +134\, q+760\, q^2+ \cdots}
\end{array}$$
where \,$\eta$\, denotes the Dedekind function on the complex upper-half plane.
The field $\Q (X_0(5))$ is generated over $\Q (X^+(5))$ by the elliptic modular function $j$, and the
relation between the functions $j$ and $t$ can be computed by using the procedure described in \cite{Go-La}\,:
$$j^2-(t^5+30\,t^4-310\,t^3-13700\,t^2-38424\,t+614000)\,j+(t^2+260\,t+5380)^3=0.$$
The $j$-invariants of quadratic $\Q$-curves of degree $5$ are obtained from this equation by rational 
values of the function $t$.

Through the above isomorphism $\Phi$, the functions \,$G(3\tau)$\, and \,$t(3\tau)$\,
can be viewed inside $\Q(X(5,3))$, as the following diagram shows:
$$\xymatrix{
  \Q\left(X(5,3)\right) \ar@{-}[d] & \Q\left(X_0(45)\right) \ar@{->}[l]_{\ \ \ \ \Phi^{\!*}}
  \ar@{-}[d]^{\pi_1^*} & \\
  \Phi^*\!\left(\Q\left(X_0(15)\right)\right) \ar@{-}[d] & \Q\left(X_0(15)\right)
  \ar@{->}[l]_{\ \ \ \ \ \simeq} \ar@{-}[rd]^{\pi_2^*} \ar@{-}[d] & \\
  \Q\left(X_0(5)\right) \ar@{-}[d] & \Q\left(G(3\tau)\right) \ar@{->}[l]_{\ \ \ \ \ \simeq} \ar@{-}[d] &
  \Q\left(X_0(5)\right) \ar@{-}[d]^{\pi_3^*} \\
  \Q\left(X^+(5)\right) & \Q\left(t(3\tau)\right) \ar@{->}[l]_{\ \ \ \ \ \simeq} & \Q\left(X^+(5)\right)
  }$$
In order to give  an explicit description of the extension $\Q\left(X(5,3)\right)/\Q\left(X^+(5)\right)$ 
on the left column of the diagram, we begin by recalling that $\pi_1$ is in fact the projection
$X_0(45)\To X_0(45)/\langle w_9\, S\, w_9\rangle$\,. We now consider the following functions
on $X_0(45)$\,:
$$U=\frac{\,(-3+ x + 9\,y)(3 + x + 9\,y)\,}{4\,x^2}\,, \ \ \ 
V= \frac{\,9(3 + x^2 + 18\,x\,y + 81\,y^2)\,}{4\,x^3}
$$
They are invariant by \,$w_9\,S\,w_9$\, and satisfy \,$[\,\Q\left(X_0(45)\right)\colon\Q\left(U,V\right)\,]=3$, so 
they generate the function field of $X_0(15)$ over $\Q$. Using the
$q$-expansions of the above functions $u$ and $v$, we obtain the following identities:
$$u= \frac{Q(U,V)}{\,2(10 + 2\,U + 3\,U^2 - 2\,V)^2\,}\,, \ \ \ \ \ 
v=\frac{R(U,V)}{\,2(10 + 2\,U + 3\,U^2 - 2\,V)^3\,}\,,$$
where
$$
\begin{array}{cl}
Q(U,V)=&-1300 -520\,U -477\,U^2 + 19\,U^3 -17\,U^4 + (260 + 52\,U + 33 \,U^2)V ,\\[5pt]
R(U,V)=&9\left(1000 - 1900\,U - 630\,U^2 - 1237\,U^3 - 39\,U^4 - 121\,U^5 + 2\,U^6\right) + \\[3pt] 
&9(- 200 + 420\,U+ 82\,U^2 + 131\,U^3 + 2\,U^4) V.
\end{array}
$$

By applying $\Phi^*$ to the resulting expressions of $u$ and $v$ in terms of $x,y$, that is, 
by changing $x,y,u,v$ \,by\, $\overline x,\overline y,\overline u,\overline v$, respectively,
in the above relations, what we get is a description 
of the subextension \,$\Q\left(X(5,3)\right)/\Phi^*\left(\Q\left(X_0(15)\right)\right)$. All we have
to do then is to give~$t$ as a rational function in $\overline u$\, and $\overline v$.
As $t(\tau)=\overline{t(3\tau)}$, this is equivalent to giving $t(3\tau)$ as a rational function 
in $u$ and $v$. Now, $G(3\tau)$ has exactly two poles at the cusps $1/5$ and $\infty$ of $X_0(15)$
with multiplicities $1$ and~$3$, respectively. Since the function
$$H(\tau) \,=\, \frac{\,\eta (3\,\tau)\,\eta (5\,\tau)^5\,}
{\,\eta(\tau)\,\eta (15\,\tau)^5\,} \,=\, \frac{1}{\,q^2\,}+\frac{1}{\,q\,}+2+2\,q+4\,q^2+\cdots$$
lies in $\Q(X_0(15))$ and has divisor $2\,(1/5) -2\,(\infty)$, we get $H=u+1$. Then, the function
$(u+1)\,G(3\tau)$ has a unique pole at $\infty$ (with multiplicity $5$), so it must be a polynomial in $u$ and
$v$. Using again the $q$-expansions of $u$ and $v$, we obtain
$$ G(3\tau) \,=\, \frac{\,u\,v-u^2-9\,u-8\,}{u+1}\,,$$
hence
\begin{equation}\label{t}
t=
\frac{\,189+205 \,\overline u+7\,\overline{u}^2+\overline{u}^3+\overline{u}^4-16\,\overline{u}\,\overline{v}-3\,\overline{u}^2\,\overline{v}\,}
{\,\overline u\,\overline v-\overline {u}^2-9\,\overline u-8\,}\,.
\end{equation}


\section{A plane quartic model for $X(5,3)_\varrho$}

The background strategy in this section is the same as in Subsection 3.1\,: a plane quartic model 
for $X(5,3)_\varrho$ over \,$\Q$\, can be theoretically obtained from a basis of the $3$-dimensional 
$\Q$-vector space $\Omega^1_\Q\!\left(X(5,3)_\varrho\right)$. So our problem amounts to
giving an explicit enough description of this space.

We recall that $\varrho$\, stands for any fixed representation of $\G_\Q$ onto $\PGL_2(\F_3)$ with
non-cyclotomic determinant. It is determined by its splitting field $L$ up to conjugation in $\PGL_2(\F_3)$,
and the condition on the determinant amounts to saying that $L$ does not contain $\sqrt{-3\,}$.
Put $K=L(\sqrt{-3\,}\,)$ and denote by $\nu$ the non-trivial element in $\Gal(K/L)$.
Since $\PGL_2(\F_3)$ is isomorphic to the symmetric group $\mathcal{S}_4$, we can take as input data
a quartic polynomial $f\in \Z[X]$ with splitting field $L$ and identify $\Gal(L/\Q)$ \,with\, $\mathcal{S}_4$
by fixing an order of the roots of~$f$. For convenience, we take as generators for this Galois group the
following permutations:\\[-15pt]
$$\sigma_1= (1,2,3), \ \ \ \ \sigma_2= (1,2) (3,4), \ \ \ \ \sigma_3=(1,2).$$\\[-20pt]

Consider on \,$\Omega^1\!\left(X(5,3)\right)=\,\Omega^1_\Q\!\left(X(5,3)\right)\otimes\Qbar$\,
the Galois action twisted by the $1$-cocycle\, $\xi$ obtained from $\varrho$\, as in Section~2.
It is defined by
$$(\overline{\omega}\otimes \gamma) ^{\sigma}_{\xi} := \left({^\sigma}\overline{\omega}\,\,\xi_\sigma^{-1}\right)
\otimes \sigma(\gamma)$$
for $\overline\omega\in\Omega^1_\Q\!\left(X(5,3)\right)$, $\gamma\in\Qbar$ and $\sigma\in\G_\Q$.
This action factors through $\Gal(K/\Q)$, and the regular differentials on $X(5,3)_\varrho$ defined over $\Q$
can be identified via the isomorphism \,$\Psi\colon X(5,3)_\varrho\To X(5,3)$\, with the fixed elements in \,$\Omega^1_K\!\left(X(5,3)\right)$\,:
$$\Omega^1_\Q\!\left(X(5,3)_\varrho\right) \,=\,
\left(\Omega^1_{\Q(\sqrt{-3\,}\,)}\!\left(X(5,3)\right)\otimes L\right)^{\Gal(K/\Q)}_\xi.$$
Moreover, the twisted action of \,$\Gal(K/\Q)$ can be restricted to the $6$-dimensional $\Q$-vector
space \,$\Omega^1_{\Q(\sqrt{-3\,}\,)}\!\left(X(5,3)\right)$, for which we take the basis
$$\left\{\overline\omega_1\,, \overline\omega_2\,, \overline\omega_3\,, \sqrt{-3\,}\,\,\overline\omega_1\,,
\sqrt{-3\,}\,\,\overline\omega_2\,, \sqrt{-3\,}\,\,\overline\omega_3 \right\}.$$
Recall that $\omega_1, \omega_2, \omega_3$ are the forms in $\Omega^1_\Q\!\left(X_0(45)\right)$ introduced
in the previous section. The action of the Galois generators \,$\sigma_1, \sigma_2, \sigma_3$\,
and \,$\nu$\, on this basis is given by the matrices
$$\begin{array}{lcl}
\tiny{s_1=\left(
\begin{matrix}
 -1/2 & 0 & \,\,0 & 0 & 0 & -3/2\\[3pt]
 -3/2  & 1 & \,\,0 & 0 & 0 & -3/2\\[3pt]
 \,\,0 & 0 & -1/2  & -3/2 & 0 & 0\\[3pt]
 \,\,0 & 0 & \,\,\,\,1/2 & -1/2 & 0 & 0\\[3pt]
 \,\,0 & 0 & \,\,\,\,1/2 & -3/2 & 1 & 0\\[3pt]
 \,1/2 & 0 & \,\,0 & 0 & 0 & -1/2\\
\end{matrix}\right),}
& \ \ &
{\tiny s_2=\left(
\begin{matrix}
    0 & 1/3 & 0 & 0 & 0 & 0 \\[3pt]
    3 & 0 & 0 & 0 & 0 & 0 \\[3pt]
    0 & 0 & -1 & 0 & 0 & 0 \\[3pt]
    0 & 0 & 0 & 0 & 1/3 & 0 \\[3pt]
    0 & 0 & 0 & 3 & 0 & 0 \\[3pt]
    0 & 0 & 0 & 0 & 0 & -1\\
\end{matrix}\right),}
\\[30pt]
{\tiny  s_3 = \begin{pmatrix}
 \phantom{c}   -1 & 0 & \phantom{c}0 &\phantom{ccc} 0 & 0 & 0 \phantom{cc}\\[3pt]
  \phantom{c}  0 & -1 &\phantom{c} 0 &\phantom{ccc} 0 & 0 & 0 \phantom{cc}\\[3pt]
   \phantom{c} 0 & 0 & \phantom{c}1 &\phantom{ccc} 0 & 0 & 0 \phantom{cc}\\[3pt]
    \phantom{c}0 & 0 &\phantom{c} 0 &\phantom{cc} -1 & 0 & 0 \phantom{cc}\\[3pt]
   \phantom{c} 0 & 0 & \phantom{c}0 &\phantom{ccc} 0 & -1 & 0\phantom{cc} \\[3pt]
   \phantom{c} 0 & 0 & \phantom{c}0 &\phantom{ccc} 0 & 0 & 1 \phantom{cc}\\
\end{pmatrix},}
& \ \ &
{\tiny s_4= \begin{pmatrix}
    -1 & 0 & 0 & 0 & \phantom{cc}0 & 0 \\[3pt]
    0 & -1 & 0 & 0 & \phantom{cc}0 & 0 \\[3pt]
    0 & 0 & 1 & 0 & \phantom{cc}0 & 0 \\[3pt]
    0 & 0 & 0 & 1 & \phantom{cc}0 & 0 \\[3pt]
    0 & 0 & 0 & 0 &\phantom{cc} 1 & 0 \\[3pt]
    0 & 0 & 0 & 0 &\phantom{cc} 0 & -1 \\
\end{pmatrix},}
\end{array}
$$
respectively. This comes from the table in Subsection~3.2 and from the definition of the
$1$-cocyle $\xi$. Indeed, the (conjugacy class of the) representation $\varrho$\, translates
into the isomorphism $\Gal(L/\Q)\simeq\Aut\left(X(5,3)\right)$ sending
\,$\sigma_1, \sigma_2, \sigma_3$ \,to\,
$\overline S^{-1}\!, \overline w_9, \overline w_5$, respectively.

We must now look for three elements in $\Omega^1_{\Q(\sqrt{-3\,}\,)}\!\left(X(5,3)\right)\otimes L$
which are linearly independent over $\Q$ and invariant by the above Galois action. To that end,
we follow the next steps:

%

\vskip 3truemm

\noindent {\bf First step.} Compute a basis for the ring of integers ${\mathcal O}_L$.

\vskip 3truemm

\noindent {\bf Second step.} Compute the $24\times 24$ matrices $\Sigma_1, \Sigma_2, \Sigma_3$ with entries in
$\Z$ giving the  action of $\sigma_1, \sigma_2, \sigma_3$, respectively, on the integral basis of ${\mathcal
O}_L$.

\vskip 3truemm

\noindent {\bf Third step.}
Form the Kronecker products
$$
\mathcal{W}_1 := s_1 \otimes\Sigma_1 \,, \quad
\mathcal{W}_2 := s_2 \otimes\Sigma_2 \,, \quad
\mathcal{W}_3 := s_3 \otimes\Sigma_3 \,, \quad
\mathcal{W}_4 := s_4 \otimes\Sigma_4 \,,
$$
where $\Sigma_4$ stands for the identity matrix $\operatorname{Id}_{24}$.
Then, compute a basis \,${X_\varrho\,,Y_\varrho\,,Z_\varrho}$\, for the $3$-dimensional vector subspace
of \,$\Omega^1_{\Q(\sqrt{-3\,}\,)}\!\left(X(5,3)\right)\otimes L$\, corresponding to the subspace
\,${\displaystyle \bigcap_{i=1}^{4}} \ker \left( \mathcal{W}_i  - \operatorname{Id}_{144}\right)$\,
of \,$\Q^{144}$.

\vskip 3truemm

\noindent {\bf Fourth step.}
Compute the $3\times 3$ matrix $\Theta$ with entries in $K $ giving the basis
change
$$(\overline{\omega}_1, \overline{\omega}_2, \overline{\omega}_3) = (X_\varrho\,,Y_\varrho\,,Z_\varrho)\,\Theta \,.$$
Replacing $X_\varrho\,, Y_\varrho\,, Z_\varrho$ by projective variables $X, Y, Z$, plugging then
$\overline{\omega}_1, \overline{\omega}_2, \overline{\omega}_3$ in
the homogenization of equation (\ref{X045}), and finally factoring out, one gets
a plane quartic equation \,$F(X,Y,Z)=0$\, for the twist $X(5,3)_\varrho$ over $\Q$.

\vskip 4truemm

\begin{rem}
The interest of working with an integral basis of ${\mathcal O}_L$ instead of just using the 
power-basis attached to a primitive element of the extension \,$L/\Q$\, is due to
the fact that one obtains experimentally better models, in the
sense of shrinking the set of bad reduction primes for the twisted
curve.
\end{rem}

\begin{rem}\label{definicio}
It can be easily checked that 
$\{ \sigma\in \Gal (K/\Q) \ | \ {}^{\sigma}\Psi=\Psi\}=\langle \nu\,\sigma_3 \rangle$. 
Thus, the isomorphism $\Psi$ is defined over $K^{\langle \nu\,\sigma_3\rangle}$ and the entries of the above
matrix~$\Theta$ belong all to this extension of degree $24$. We also have 
$$\{\sigma\in\Gal (K/\Q) \ | \ \,{}^{\sigma}\Psi\,
\Psi^{-1}\in\langle\overline w_9\,\overline S\,\overline
w_9\,\rangle\,\} \,=\, \langle \nu\,\sigma_3, \,\sigma_2\,\sigma_1\,\sigma_2 \rangle\,,$$ 
which means that the composition
$$
X(5,3)_\rho\,\stackrel{\Psi}\To X(5,3)\,\To\, X(5,3)/
{\langle\overline w_9\,\overline S\,\overline w_9 \rangle}\,\simeq\,X_0(15)$$ 
is defined over a number field of degree $8$ which is the compositum of the quadratic field $k$ 
attached to the character \,$\eps\det\varrho$\, and the number field $L_0$ generated
by the root of the quartic polynomial $f$ fixed by the permutation \,$\sigma_2\,\sigma_1\,\sigma_2$.
\end{rem}


\section{An explicit example}

All the steps in the previous section can be performed using a software package for algebraic 
manipulation such as Pari or Magma. As an example, we consider the surjective Galois representation
$\varrho\,\colon\G_\Q\to\PGL_2(\F_3)$ defined up to conjugation by the splitting field of
the irreducible polynomial
$$f(X)=X^4-3\,X^2 + 2\,X +3\,.$$
Since the discriminant of $f$ equals $-33$ up to squares, the
field $k$ over which the quadratic $\Q$-curves realizing
$\varrho$ must be defined is $\Q(\sqrt{11}\,)$  (see Section~2).

Note that if the rank of $X_0(15)(k\,L_0)$ is zero, where $L_0$ is
the quartic number field defined by $f$ (cf. Remark~\ref{definicio}), then we do not
need to compute an equation for $X(5,3)_\varrho$. Indeed,  in this
case the values $t\in X^+(5)(\Q)$ obtained from the torsion points
of $X_0(15)(k\,L_0)$ using (\ref{t}) would provide  us with  a
finite set of candidate $\Q$-curves $E$, and then it would suffice to
check whether \,$\varrho_{E,3}=\varrho$\, or not for each of them. In our example, this 
strategy does not apply, since the rank of $X_0(15)(k)$, which can be computed using 
Magma V2.11, is one.

The procedure described in the previous section allows us to
obtain the following projective model for $X(5,3)_\varrho$\,:
$$\begin{array}{l}
 -9\,X\,Y\,( 2\,X + Y )\,( 9\,X + 8\,Y )
 + 9\,( 6\,X^3 + 62\,X^2\,Y + 66\,X\,Y^2 + 15\,Y^3)\, Z \,+ \\[8pt] 
 \,\,\,\,3\,( 27\,X^2 - 104\,X\,Y - 83\,Y^2 ) \, Z^2 -3\,( 94\,X + 7\,Y) \, Z^3  + 191\, Z^4 \,=\, 0\,.
\end{array}$$
We have found four rational points on the line at infinity, namely
$$
P_1= [0 \colon  1 \colon 0]  \,, \quad P_2= [1 \colon 0 \colon 0]
\,, \quad P_3= [1 \colon -2 \colon 0] \,, \quad P_4= [8 \colon
-9\colon  0] \,.
$$
The $j$-invariants, up to Galois conjugation, for the corresponding $\Q$-curves of 
degree~$5$ realizing $\varrho$ are
$$
\begin{array}{ccl}
j_1&=&\left(-8 \sqrt{11}\,\right)^3(10 + 3 \sqrt{11}\,)\,,   \\ [6pt]
j_2&=& \left(6(110 + 31 \sqrt{11}\,)\right)^3(10 + 3 \sqrt{11}\,)\,,  \\[6pt]
j_3&=& \left(12 (10 - \sqrt{11}\,)\right)^3(10 + 3 \sqrt{11}\,)\,,  \\[6pt]
j_4&=&\left(2 (-6878815950 +2118474913\sqrt{11}\,)/53^5\right)^3 (10 + 3\sqrt{11}\,)\,.\\[6pt]
\end{array}
$$

Let us finish by showing that the Chabauty-Coleman method to determine the set of rational points 
on a curve of genus at least two cannot be applied to~$X(5,3)_\varrho$\, for the representation $\varrho$\,
in this example, since the requirement that the rank of the jacobian $J\left(X(5,3)_\varrho\right)$ be 
smaller than the genus of the curve is not fulfilled. Indeed, consider the projection 
\,$X(5,3)_\varrho\To X_0(15)$\, in Remark \ref{definicio} and the corresponding images $Q_1, Q_2, Q_3, Q_4$ 
in $X_0(15)(k\,L_0 )$ of the above points $P_1, P_2, P_3, P_4$. The morphism
$$
X(5,3)_\varrho\,\hookrightarrow\, J\left(X(5,3)_\varrho\right)\,,\qquad P\,\mapsto\, (P)-(P_1)
$$
is defined over $\Q$ and produces three rational points in
$J\left(X(5,3)_\varrho\right)(\Q)$ whose images in the elliptic curve $J_0(15)$ are $Q_2-Q_1$, $Q_3-Q_1$ 
and  $Q_4-Q_1$. 
Now, using Magma again, one can compute the $3\times 3$
matrix obtained by applying the N\'eron-Tate pairing 
$$<P,Q> \,=\, \frac{1}{\,2\,}\left(h(P+Q)-h(P)-h(Q)\right)$$ 
over the points $Q_2-Q_1$, $Q_3-Q_1$, $Q_4-Q_1$. An approximation for the determinant of this 
matrix turns out to be $6.460235$. Since it is nonzero, these three points are linearly
independent on $J_0(15)$. It follows that the rank of 
$J\left(X(5,3)_\varrho\right)(\Q)$ is at least three.

\bibliographystyle{alpha}
\bibliography{pqt}

\begin{thebibliography}{FLR02}

\bibitem[BGGP]{bggp}
M.~Baker, E.~Gonz\'alez, J.~Gonz\'alez, and B.~Poonen.
\newblock Finiteness results for modular curves of genus at least $2$.
\newblock To appear in \emph{Amer. J. Math.}

\bibitem[Cre97]{cremona97}
J.~E. Cremona.
\newblock {\em Algorithms for modular elliptic curves}.
\newblock Cambridge University Press, 1997.

\bibitem[Fer03]{julio}
J.~Fern{\'a}ndez.
\newblock {\em Elliptic realization of {G}alois representations}.
\newblock PhD thesis, Universitat Polit{\`e}cnica de Catalunya, 2003.

\bibitem[Fer04]{Fe}
J.~Fern{\'a}ndez.
\newblock A moduli approach to quadratic $\mathbb{Q}$-curves realizing
  projective mod $p$\,\, {G}alois representations.
\newblock 2004.
\newblock Preprint available at \texttt{http://www.math.leidenuniv.nl/gtem}.

\bibitem[FLR02]{Fe-La-Rio}
J.~Fern{\'a}ndez, J.-C. Lario, and A.~Rio.
\newblock Octahedral {G}alois representations arising from $\mathbb {Q}$-curves
  of degree $2$.
\newblock {\em Canad. J. Math.}, 54(6):1202--1228, 2002.

\bibitem[GL98]{Go-La}
J.~Gonz{\'a}lez and J.-C. Lario.
\newblock Rational and elliptic parametrizations of $\mathbb {Q}$-curves.
\newblock {\em J. Number Theory}, 72(1):13--31, 1998.

\bibitem[KM88]{Ke-Mo}
M.~A. Kenku and F.~Momose.
\newblock Automorphism groups of the modular curves {$X\sb 0(N)$}.
\newblock {\em Compositio Math.}, 65(1):51--80, 1988.

\bibitem[Lig77]{Li}
G.~Ligozat.
\newblock Courbes modulaires de niveau {$11$}.
\newblock In {\em Modular functions of one variable V}, pages 149--237. {\rm
  Lecture Notes in Math., vol. 601}. Springer, 1977.

\bibitem[LN64]{Le-Ne}
J.~Lehner and M.~Newman.
\newblock Weierstrass points on {$\Gamma\sb 0(n)$}.
\newblock {\em Ann. of Math.}, 79:360--368, 1964.

\bibitem[Maz77]{Ma77}
B.~Mazur.
\newblock Rational points on modular curves.
\newblock In {\em Modular functions of one variable V}, pages 107--148. {\rm
  Lecture Notes in Math., vol. 601}. Springer, 1977.

\bibitem[Ogg74]{Ogg}
A.~P. Ogg.
\newblock Hyperelliptic modular curves.
\newblock {\em Bull. Soc. Math. France}, 102:449--462, 1974.

\end{thebibliography}

\vfill

\begin{tabular}{l@{\ \hskip 1 truecm\ }l@{\ \hskip 1 truecm\ }l}
Julio Fern\'andez & Josep Gonz\'alez & Joan-C. Lario \\
julio@mat.upc.edu  & josepg@mat.upc.edu & joan.carles.lario@upc.edu
\end{tabular}

\vskip 0.35truecm

\begin{tabular}{l}
Facultat de Matem\`atiques i Estad\'\i stica  \\
Universitat Polit\`ecnica de Catalunya  \\
Pau Gargallo 5  \\
08028 Barcelona, Spain
\end{tabular}

\end{document}